\documentclass[11pt, reqno]{amsart}
\usepackage{amsmath, amsthm, amscd, amsfonts, amssymb, graphicx, color}
\usepackage[bookmarksnumbered, colorlinks, plainpages]{hyperref}

\newtheorem{theorem}{Theorem}[section]
\newtheorem{lemma}[theorem]{Lemma}

\theoremstyle{definition}

\theoremstyle{remark}

\numberwithin{equation}{section}

\input{mathrsfs.sty}

\begin{document}
\setcounter{page}{1}

\title{Commutators of automorphic  composition operators with adjoints}

\author{Liangying Jiang}

\address{Department of Applied Mathematics,  Shanghai Finance University, Shanghai 201209, P. R. China}

\email{\textcolor[rgb]{0.00,0.00,0.84}{liangying1231@163.com,
jiangly@shfc.edu.cn}}



\date{ May 12th, 2014
\newline \indent Supported by the National
Natural Science Foundation of China (No.11101279)}

\begin{abstract} In this paper, we investigate the compactness of the commutator $[C_\psi^{\ast}, C_\varphi]$  on the Hardy space $H^2(B_N)$
or the weighted Bergman space  $A^2_s(B_N)$ ($s>-1$), when $\varphi$ and $\psi$ are automorphisms of the unit ball $B_N$. We obtain that $[C_\psi^{\ast}, C_\varphi]$
is compact if and only if $\varphi$ and $\psi$ commute and they are both unitary. This generalizes the corresponding result in one variable. Moreover,
our technique is different and  simpler. In addition, we also discuss the commutator $[C_\psi^{\ast}, C_\varphi]$  on the Dirichlet space $\mathcal{D}(B_N)$, where $\varphi$ and $\psi$ are linear fractional self-maps  or both  automorphisms of $B_N$.
\end{abstract}
\maketitle

\section{Introduction}

Let $B_N$ denote the unit ball of $\mathbb{C}^N$ and let $\varphi$ be a holomorphic self-map of $B_N$. We define the composition operator $C_\varphi$
 by $C_\varphi(f)=f\circ \varphi$, where $f$ is analytic in $B_N$.

In this paper, we are interested in characterizing the compactness of the  commutator $[C_\psi^{\ast}, C_\varphi]=C_\psi^{\ast} C_\varphi- C_\varphi C_\psi^{\ast}$ on some classical function spaces, when  $\varphi$ and $\psi$ are automorphisms of $B_N$. The motivation comes from  a recent work of
Clifford et al. \cite{CLN12}. They considered when the  commutator $[C_\psi^{\ast}, C_\varphi]$ is non-trivially compact on the Hardy space $H^2(D)$
for linear fractional self-maps  $\varphi$ and $\psi$ of the unit disk $D$. Here, non-trivially compact means that $[C_\psi^{\ast}, C_\varphi]$ is compact
but nonzero, moreover, both $C_\psi^{\ast} C_\varphi$ and $C_\varphi C_\psi^{\ast}$ are not compact. In particular, when  $\varphi$ and $\psi$ are
automorphisms of $D$, they showed that  $[C_\psi^{\ast}, C_\varphi]$ is non-trivially compact if and only if both maps are rotations. All results
were extended by MacCluer et al. \cite{MNW13} to the weighted Bergman space $A^2_s(D)$ ($s>-1$).

In Section 3, we first investigate the commutator $[C_\psi^{\ast}, C_\varphi]$  on the Hardy space $H^2(B_N)$
  and the weighted Bergman space  $A^2_s(B_N)$ ($s>-1$), where $\varphi$ and $\psi$ are automorphisms of $B_N$, neither  $\varphi$ nor $\psi$ is the idetity.
In this case, we will prove that  $[C_\psi^{\ast}, C_\varphi]$ is compact if and only if  $\varphi$ and $\psi$ commute and  both maps are unitary. This
generalizes the result of  automorphisms case in the unit disk to the unit ball. In order to deduce that
 $\varphi$ and $\psi$ commute when $[C_\psi^{\ast}, C_\varphi]$ is compact on  $A^2_s(D)$, MacCluer et al. \cite{MNW13} had done lots of  complicated
calculations. It is difficult for us to use this similar idea. So in Section 3, we will find another simpler technique to solve similar problem in higher
dimensions, which involved in the application of the semi-multiplication property for Toeplitz operators.

Furthermore, we will extend to discuss the compactness of  $[C_\psi^{\ast}, C_\varphi]$ on  the Dirichlet space $\mathcal{D}(B_N)$ in Section 4.
Based on the adjoint formula for $C_\psi$ on $\mathcal{D}(B_N)$, we   obtain the necessary and sufficient condition for $[C_\psi^{\ast}, C_\varphi]$ to be compact in terms of linear fractional self-maps $\varphi$ and $\psi$ of $B_N$.  As an immediate result, for automorphisms $\varphi$ and $\psi$ of $B_N$,  we find that   $[C_\psi^{\ast}, C_\varphi]$ is compact
on  $\mathcal{D}(B_N)$  if and only if $\varphi$ and $\psi$ commute.  More specially, when
 $\varphi$ and $\psi$ are linear fractional self-maps of $D$, the condition for $[C_\psi^{\ast}, C_\varphi]$ to be compact on the  Dirichlet space $\mathcal{D}(D)$ is the same as that on the Hardy space  $H^2(D)$ and  the weighted Bergman space $A^2_s(D)$ ($s>-1$).

\section{ Preliminaries}

Here we collect some necessary background information.

\subsection{Analytic function spaces}
Let $\partial B_N$ denote the boundary of the unit ball $B_N$. The Hardy space $H^2(B_N)$ is defined by
$$H^2(B_N)=\{f\ \mbox{analytic in}\ B_N: ||f||^2\equiv \sup\limits_{0<r<1}\int_{\partial B_N}|f(r\zeta)|^2
d\sigma(\zeta)<\infty\},$$ where $d\sigma$ denotes the normalized surface measure on $\partial B_N$. The weighted Bergman space
$A^2_s(B_N)$, for $s>-1$, is defined by
$$A^2_s(B_N)=\{f\ \mbox{analytic in}\ B_N: ||f||_s^2\equiv \int_{B_N}|f(z)|^2 dv_s(z)<\infty\},$$
where $dv_s(z)=\frac{\Gamma(N+s+1)}{N!\Gamma(s+1)}(1-|z|^2)^sdv(z)$ and  $dv$ denotes the normalized volume measure on $B_N$. In this paper,
we will often use  $\mathcal{H}$ to denote the Hardy space $H^2(B_N)$ or the weighted Bergman space $A^2_s(B_N)$. It is well known that both
the Hardy space and  the weighted Bergman space are the reproducing kernel Hilbert spaces, where the reproducing kernel is given by
$K_w(z)=(1-<z, w>)^{-t}$ for $z, w \in B_N$, with $t=N$ for $H^2(B_N)$  and $t=N+s+1$ for $A^2_s(B_N)$. So the the normalized reproducing kernel
is given by  $$k_w(z)=\frac{K_w(z)}{||K_w||_{\mathcal{H}}}=\frac{(1-|w|^2)^{t/2}}{(1-<z, w>)^t}.$$

A multi-index $\alpha=(\alpha_1,\ldots,\alpha_N)$ is an $N$-tuple of non-negative integers $\alpha_i$. The total order of a multi-index
 is given by $|\alpha|=|\alpha_1|+\cdots+|\alpha_N|$. Let $\alpha!=\alpha_1!\cdots \alpha_N!$ and $z^\alpha=z_1^{\alpha_1}\cdots z_N^{\alpha_N}$
 for $z=(z_1,\ldots,z_N)\in \mathbb{C}^N$. An analytic function $f$ in $B_N$ has a power series representation $$f(z)=\sum\limits_\alpha c_\alpha z^\alpha,$$
where the sum is over all multi-indexes.

The Dirichlet space $\mathcal{D}(B_N)$ is defined as
$$\mathcal{D}(B_N)=\{f(z)=\sum\limits_\alpha c_\alpha z^\alpha\, \mbox{analytic in}\, B_N: ||f||_{\mathcal{D}^0}^2\equiv \sum\limits_\alpha |c_\alpha |^2|\alpha|\frac{\alpha!}{|\alpha|!}<\infty\},$$
where the quantity $||\cdot||_{\mathcal{D}^0}$ defines a semi-norm on $\mathcal{D}(B_N)$. We equip it with the norm $$||f||_{\mathcal{D}}^2=|f(0)|^2+||f||_{\mathcal{D}^0}^2$$
and the inner product $<\cdot , \cdot>_{\mathcal{D}}$. So the reproducing kernel for  $\mathcal{D}(B_N)$ is given by
 $$K_w(z)=1+\log \frac{1}{1-<z,w>},\quad z, w\in B_N.$$

When acting on the reproducing kernel, composition operators and Toeplitz operators have the following adjoint property:
$$C_\varphi^\ast T_h^\ast K_w=\overline{h(w)}K_{\varphi(w)}, \quad w\in B_N$$ for all analytic self-maps $\varphi$ of $B_N$ and $h\in L^\infty(B_N)$.

\subsection{Adjoint formula}

Given a bounded measurable complex-valued function $b$ on $\partial B_N$ (or  $B_N$), the Toeplitz operator $T_b$ on $H^2(B_N)$ (or $A^2_s(B_N)$) is defined by $$T_b f=P(bf),$$ where $P$ is the orthogonal projection from $L^2(\partial B_N)$ (or $L^2(B_N, dv_s)$) onto $H^2(B_N)$ (or $A^2_s(B_N)$). If $b$ is analytic, then $T_b$ is a multiplication by $b$.

In this paper, we will often  use the  semi-multiplication property for Toeplitz operators mod $\mathcal{K}$,  where $\mathcal{K}$ denotes the ideal  of compact operators.  That is, if $b\in L^\infty(\partial B_N)$ (or $ L^2(B_N)$) and $h\in C(\overline{B_N})$, then $T_bT_h-T_{bh}$ is compact on $H^2(B_N)$ (or $A^2_s(B_N)$). This result on the Hardy space $H^2(B_N)$ is Proposition 1.4 in \cite{McD77}; the Bergman space version comes from the proof of Theorem 1 in \cite{Cob73}; similar to the proof of Theorem 1 in \cite{Cob73}, this property also can be extended to the weighted Bergman space  $A^2_s(B_N)$ ($s>-1$). In the unit disk $D$, this fact has been well described in page 73 of \cite{MNW13}.
\\ \\
{\bf Theorem A.}
{\it  Suppose that $$\varphi(z)=\frac{az+b}{cz+d}$$ is a linear fractional self-map of $D$ and $ad-bc\ne 0$. Then the adjoint of $C_\varphi$ acting on the Hardy space $H^2(D)$ or the weighted Bergman space $A^2_s(D)$ ($s>-1$) is given by $$C_{\varphi}^{\ast}=T_g C_{\sigma} T_h^{\ast},$$
where  $$\sigma(z)=\frac{\overline{a}z-\overline{c}}{-\overline{b}z+\overline{d}}$$  is the Krein adjoint of $\varphi$,  $g(z)=(-\overline{b}z+ \overline{d})^{-t}$  and
$h(z)=(cz+ d)^t$ are in $H^\infty$ with $t=1$ for $H^2(D)$  and $t=s+2$  for $A^2_s(D)$.}
\\ \par
This adjoint formula of $C_\varphi$ was first established on $H^2(D)$  by Cowen  \cite{Co88}  and was generalized to
$A^2_s(D)$ ($s>-1$) by Hurst \cite{Hu97}. In the unit ball $B_N$, Cowen and MacCluer \cite{CM00} obtained the following similar  adjoint formula for  $C_\varphi$ when acting on $H^2(B_N)$  or   $A^2_s(B_N)$ ($s>-1$).
\\ \\
{\bf Theorem B.}
{\it Let $$\varphi(z)=\frac{Az+B}{<z, C>+d}$$
be a linear fractional self-map of $B_N$, where $A$ is an $N\times N$ matrix, $B$ and $C$
  are $N\times 1$ matrices,  and $d$ is a scalar. Then on the space $\mathcal{H}$,  $$C_{\varphi}^{\ast}=T_g C_{\sigma} T_h^{\ast},$$
where $$\sigma(z)=\frac{A^\ast z-C}{<z,  -B>+ \bar{d}},$$ $g(z)=(<z,  -B>+ \bar{d})^{-t}$  and
$h(z)=(<z, C>+ d)^t$ with $t=N$ when $\mathcal{H}=H^2(B_N)$  and $t=N+s+1$
when  $\mathcal{H}=A^2_s(B_N)$ ($s>-1$).}
\\ \par
We will refer to the functions $g, h$ and $\sigma$ as the auxiliary functions of $\varphi$ when they connected by the equation $C_{\varphi}^{\ast}=T_g C_{\sigma} T_h^{\ast}$. We also need frequently use the  property: if $\varphi\in\mbox{Aut}(B_N)$ (or $\mbox{Aut}(D)$) then $\sigma=\varphi^{-1}$,
where $\mbox{Aut}(B_N)$ (or $\mbox{Aut}(D)$) denotes the set of automorphisms of $B_N$ (or $D$).

Moreover, we have a different adjoint formula for $C_{\varphi}$ on the Dirichlet space $\mathcal{D}(B_N)$ as the following.
\\ \\
{\bf Theorem C.} (Theorem 7 in \cite{Po11})
{\it Let $\varphi$ be a linear fractional self-map of $B_N$ and $K_w$ be the reproducing  kernel  of $\mathcal{D}(B_N)$. Then for $f\in \mathcal{D}(B_N)$, we have
 $$C_{\varphi}^{\ast}f=f(0)K_{\varphi(0)}+C_\sigma f-f(\sigma(0)),$$ where
$\sigma$ is the Krein adjoint of $\varphi$.}
\\ \par This means that the adjoint  of $C_{\varphi}$ on  $\mathcal{D}(B_N)$ can be identified as another composition operator and a rank $2$ operator. For the unit disk case, this adjoint formula on the Dirichlet space $\mathcal{D}(D)$ was given by Gallardo-Guti\'{e}rrez and  Montes-Rodr\'{i}guez (see Theorem 3.3 in \cite{GM03}).

It is well known that when $\varphi$ is a linear fractional self-map of $B_N$, then $C_\varphi$ is compact on the space mentioned in this paper if
and only if $||\varphi||_\infty<1$. Let $\varphi$ be  a linear fractional self-map of $D$ with a fixed point $\omega\in\partial D$, if $\omega$
is the only fixed point for $\varphi$, we say that $\varphi$ is parabolic. If
$\varphi$ has an additional  fixed point, we call $\varphi$  hyperbolic.

\subsection{Julia-Carath\'{e}odory theorem}
Given $\zeta\in \partial B_N$, a continuous function $\gamma:[0, 1]\rightarrow B_N$ with $\lim\limits_{t\to 1}\gamma(t)=\zeta$
is said to a restricted $\zeta$-curve if
$$\lim\limits_{t\to 1} \frac{|\gamma(t)-<\gamma(t),\zeta>|^2}{1-|<\gamma(t),\zeta>|^2}=0
\quad  \mbox{and} \quad  \sup\limits_{0\le t<1}\frac{|\zeta-<\gamma(t),\zeta>\zeta|}{1-|<\gamma(t),\zeta>|}<\infty.$$
If $\lim\limits_{t\to 1}f(\gamma(t))=f(\zeta)$ for  every restricted $\zeta$-curve $\gamma$, we say that $f:  B_N\rightarrow \mathbb{C}$
 has restricted limit  and  write $R\lim\limits_{z\to\zeta}f(z)=f(\zeta)$.

Let $\varphi$ be a holomorphic  self-map of $B_N$. We say that $\varphi$ has finite angular
derivative at $\zeta\in \partial B_N$,
if there exists  a point  $\eta\in \partial B_N$ so that
 $$A_\varphi(\zeta)=R\lim\limits_{z\to\zeta}\frac{<1-<\varphi(z), \eta>}{ 1-<z,
\zeta>}$$  exists. We write $\varphi_\eta:=<\varphi, \eta>$ and $D_\zeta=\frac{\partial}{\partial \zeta}$ for the directional derivative in the direction of $\zeta$, and we put $$d_{\varphi}(\zeta)=\liminf\limits_{z\to
\zeta}\frac{1-|\varphi(z)|}{1-|z|}$$

The following is the Julia-Carath\'{e}odory theorem for the ball (see Theorem 8.5.6 in \cite{Rudin} or Theorem 2.2 in \cite{CKP14})
\\ \\
{\bf Theorem D.}
{\it  Let $\varphi$ be a holomorphic  self-map of $B_N$ and $\zeta\in \partial B_N$. The following statements are equivalent:
\\ (1)  $\varphi$ has finite angular
derivative at $\zeta$.
\\ (2) $d_{\varphi}(\zeta)<\infty$.
\\ (3)  $\varphi$ has restricted limit   $\eta\in \partial B_N$  at $\zeta$ and  $D_{\zeta}\varphi_{\eta}(z)=<\varphi'(z)\zeta, \eta>$  has
finite restricted

limit at $\zeta$.
\\ Moreover, when these conditions hold, the following statements hold:
\\  (4) $D_{\zeta}\varphi_{\eta}(z)=<\varphi'(z)\zeta, \eta>$  has
restricted limit at $\zeta$ with  $D_{\zeta}\varphi_{\eta}(\zeta)=d_{\varphi}(\zeta)$.
\\   (5)  $A_\varphi(\zeta)=d_{\varphi}(\zeta)$.
\\  (6)  $\frac{\varphi_{\eta^\bot}(z)}{1-<z, \zeta>}$ has restricted limit $0$ at $\zeta$ for any $\eta^\bot\in \partial B_N$ orthogonal to $\eta$.}


\section { The commutator on the Hardy space and the weighted Bergman space}

In this section, we discuss the compactness of the commutator $[C_\psi^{\ast}, C_\varphi]$  on the Hardy space $H^2(B_N)$
  and the weighted Bergman space  $A^2_s(B_N)$ ($s>-1$),  when $\varphi$ and $\psi$ are automorphisms of  $B_N$. We will show the
following main theorem.

\begin{theorem} \label{the 3.1}  Let $\varphi$  and $\psi$ be  automorphisms of  $B_N$, neither  $\varphi$  nor  $\psi$ is the identity. The commutator $[C_\psi^{\ast}, C_\varphi]$ is compact on  $H^2(B_N)$  or $A^2_s(B_N)$ ($s>-1$) if and only if  $\varphi$ and $\psi$  commute and both maps are unitary.
\end{theorem}

Before we give a proof for Theorem 3.1, we need some useful lemmas. The first lemma is analogous to Lemma 3.2 of \cite{CLN12} and Lemma 4.3 of \cite{MNW13}.
But the result has been improved. Recall that for $w\in B_N$, the function $k_w$ is the normalized reproducing kernel given by $k_w=K_w/||K_w||_{\mathcal{H}}$.

\begin{lemma}\label{lem 3.2}  Assume that $\varphi$ and $\psi$ are holomorphic self-maps of $B_N$. Suppose that there exist points $\zeta_1$ and $\zeta_2$ on $\partial B_N$ such that $\varphi(\zeta_1)=\psi(\zeta_2)=\omega\in \partial B_N$ and  $A_\varphi(\zeta_1)$ and $A_\psi(\zeta_2)$
exist. Then on the space $\mathcal{H}$,
$$\lim\limits_{r\to 1}<C_{\psi}^{\ast} k_{r\zeta_2}, C_{\varphi}^{\ast} k_{r\zeta_1}>=\biggl(\frac{2}{d_\varphi(\zeta_1)+d_\psi(\zeta_2)}\biggr)^t>0,$$
where $t=N$ when  $\mathcal{H}=H^2(B_N)$  and $t=N+s+1$
when  $\mathcal{H}=A^2_s(B_N)$ ($s>-1$).
\end{lemma}

\proof Let $U$ be a unitary map of $B_N$ so that $U\omega=e_1$. Then $U\varphi(\zeta_1)=U\psi(\zeta_2)=U\omega=e_1$. Write $\phi=U\varphi$ and $\rho=U\psi$,
we see that $$<C_{\rho}^{\ast} k_{r\zeta_2}, C_{\phi}^{\ast} k_{r\zeta_1}>=<C_{U}^{\ast}C_{\psi}^{\ast} k_{r\zeta_2}, C_{U}^{\ast}C_{\varphi}^{\ast} k_{r\zeta_1}>=<C_{\psi}^{\ast} k_{r\zeta_2}, C_{\varphi}^{\ast} k_{r\zeta_1}>.$$
Thus we may assume $\omega=e_1$. Note that $$\frac{1}{<C_{\psi}^{\ast} k_{r\zeta_2}, C_{\varphi}^{\ast} k_{r\zeta_1}>}=\frac{||K_{r\zeta_1}||_{\mathcal{H}}||K_{r\zeta_2}||_{\mathcal{H}}}{<K_{\psi(r\zeta_2)}, K_{\varphi(r\zeta_1)}>}=\biggl(\frac{1-<\varphi(r\zeta_1), \psi(r\zeta_2)>}{1-r^2}\biggr)^t$$
and{\setlength\arraycolsep{2pt}
\begin{eqnarray*}\frac{1-<\varphi(r\zeta_1), \psi(r\zeta_2)>}{1-r^2}
&=&\frac{1-|\psi(r\zeta_2)|^2}{1-r^2}+\frac{|\psi(r\zeta_2)|^2-<\varphi(r\zeta_1), \psi(r\zeta_2)>}{1-r^2}
\\ &=& \frac{1-|\psi(r\zeta_2)|^2}{1-r^2}+\frac{<\psi(r\zeta_2)-\varphi(r\zeta_1), \psi(r\zeta_2)>}{1-r^2}.
\end{eqnarray*}}Now, we calculate that{\setlength\arraycolsep{2pt}
\begin{eqnarray*} && \frac{<\psi(r\zeta_2)-\varphi(r\zeta_1), \psi(r\zeta_2)>}{1-r^2}
\\ &=& \frac{\psi_1(r\zeta_2)-\varphi_1(r\zeta_1)}{1-r^2}\cdot\overline{\psi_1(r\zeta_2)}+
\sum_{j=2}^N\frac{\psi_j(r\zeta_2)-\varphi_j(r\zeta_1)}{1-r^2}\cdot\overline{\psi_j(r\zeta_2)}
\\ &=&\biggr(\frac{1-<\varphi(r\zeta_1), e_1>}{1-<r\zeta_1, \zeta_1>} -\frac{1-<\psi(r\zeta_2), e_1>}{1-<r\zeta_2, \zeta_2>}\biggl)\frac{\overline{\psi_1(r\zeta_2)}}{1+r}
\\ && +
\sum_{j=2}^N \biggr(\frac{\psi_j(r\zeta_2)}{1-<r\zeta_2, \zeta_2>} -\frac{\varphi_j(r\zeta_1)}{1-<r\zeta_1, \zeta_1>}\biggl)\frac{\overline{\psi_j(r\zeta_2)}}{1+r}.
\end{eqnarray*}}Since $\varphi$ and $\psi$ have finite angular derivatives respectively at $\zeta_1$ and $\zeta_2$, by Theorem D(5), we get that
$$\lim\limits_{r\to 1}\frac{1-<\varphi(r\zeta_1), e_1>}{1-<r\zeta_1, \zeta_1>}=A_\varphi(\zeta_1)=d_\varphi(\zeta_1)$$
and  $$\lim\limits_{r\to 1}\frac{1-<\psi(r\zeta_2), e_1>}{1-<r\zeta_2, \zeta_2>}=A_\psi(\zeta_2)=d_\psi(\zeta_2).$$
Moreover, by Theorem D(6),
$$\lim\limits_{r\to 1}\frac{\varphi_j(r\zeta_2)}{1-<r\zeta_1, \zeta_1>}=0 \qquad
\mbox{and} \qquad \lim\limits_{r\to 1}\frac{\psi_j(r\zeta_2)}{1-<r\zeta_2, \zeta_2>}=0$$
hold for $2\le j\le N$. Therefore, {\setlength\arraycolsep{2pt}
\begin{eqnarray*}&& \lim\limits_{r\to 1}\frac{1}{<C_{\psi}^{\ast} k_{r\zeta_2}, C_{\varphi}^{\ast} k_{r\zeta_1}>}=\lim\limits_{r\to 1}\biggl(\frac{1-<\varphi(r\zeta_1), \psi(r\zeta_2)>}{1-r^2}\biggr)^t
\\ &=&  \lim\limits_{r\to 1} \biggl(\frac{1-|\psi(r\zeta_2)|^2}{1-r^2}+\frac{<\psi(r\zeta_2)-\varphi(r\zeta_1), \psi(r\zeta_2)>}{1-r^2}\biggr)^t
\\ &=&\biggr[d_\psi(\zeta_2)+\frac{1}{2}(d_\varphi(\zeta_1)-d_\psi(\zeta_2))\biggl]^t
\\ &=&\biggl(\frac{d_\varphi(\zeta_1)+d_\psi(\zeta_2)}{2}\biggr)^t.
\end{eqnarray*}}and we obtain the desire conclusion.
\ \ $\Box$

\begin{lemma}\label{lem 3.3} Assume that $\varphi$ and $\psi$ are holomorphic self-maps of $D$. Suppose that there exist points $\zeta_1$ and $\zeta_2$ on $\partial D$ such that $\varphi(\zeta_1)=\psi(\zeta_2)=\omega\in \partial D$ and the angular derivatives $\varphi'(\zeta_1)$ and $\psi'(\zeta_2)$
exist. Then on $H^2(D)$ or $A^2_s(D)$ ($s>-1$),
$$\lim\limits_{r\to 1}<C_{\psi}^{\ast} k_{r\zeta_2}, C_{\varphi}^{\ast} k_{r\zeta_1}>=\biggl(\frac{2}{|\varphi'(\zeta_1)|+|\psi'(\zeta_2)|}\biggr)^t>0,$$
where $$k_w(z)=\biggl(\frac{\sqrt{1-|w|^2}}{1-z\overline{w}}\biggr)^t,\qquad z, w\in D$$ is the normalized reproducing kernel
with  $t=1$ for $H^2(D)$  and $t=s+2$
for $A^2_s(D)$.
\end{lemma}

This is an easy corollary of Lemma 3.2. In fact, Lemma 3.3 can also be immediately obtained  if we we notice that $\zeta_1\varphi'(\zeta_1)\omega=|\varphi'(\zeta_1)|$ and $\zeta_2\psi'(\zeta_2)\omega=|\psi'(\zeta_2)|$ from the Julia-Carath\'{e}odory  Theorem of the unit disk (see Theorem 2.44 of \cite{CM1}). Because  Lemma 3.2 of \cite{CLN12} and Lemma 4.3 of \cite{MNW13} have given that
$$\lim\limits_{r\to 1}<C_{\psi}^{\ast} k_{r\zeta_2}, C_{\varphi}^{\ast} k_{r\zeta_1}>=\biggl(\frac{2\omega}{2\omega|\psi'(\zeta_1)|+[\zeta_1\varphi'(\zeta_1)-\zeta_2\psi'(\zeta_2)|}\biggr)^t,$$
we only need multiple $\overline{\omega}$ on numerator and denominator of the above fraction to deduce the limit in Lemma 3.3.

We now show that when  $\varphi$ and $\psi$ are automorphisms  of $D$, in order for the commutator $[C_\psi^{\ast}, C_\varphi]$ to be compact, the inducing
maps $\varphi$ and $\psi$ must commute. This result on the Hardy space $H^2(D)$ is Lemma 5.1 of \cite{CLN12}, and it on the weighted Bergman space
$A^2_s(D)$ is Theorem 4.6 of \cite{MNW13}. Their technique is similar, but on $A^2_s(D)$, some calculations involved are very complex. It is
almost impossible to do similar calculations in higher dimensions. We have to find other method which can be used on  $H^2(D)$ and  $A^2_s(D)$ simultaneously
and which can also be extended to the unit ball.

\begin{lemma}\label{lem 3.4} Assume that $\varphi$ and $\psi$ are automorphisms  of $D$. If the commutator $[C_\psi^{\ast}, C_\varphi]$ is compact on  $H^2(D)$ or  $A^2_s(D)$ ($s>-1$), then $\varphi$ and $\psi$ commute.
\end{lemma}

\proof This is Lemma 5.1 in \cite{CLN12} and Theorem 4.6 in \cite{MNW13}. We will give another simpler proof and we only focus on the Hardy space $H^2(D)$.
Set $$\varphi(z)=\frac{a_1z+b_1}{c_1z+d_1}\qquad \mbox{and} \qquad \psi(z)=\frac{a_2z+b_2}{c_2z+d_2}$$
with normalization $a_1d_1-c_1b_1=1$ and $a_2d_2-c_2b_2=1$.
By Theorem A, we have  $$C_{\varphi}=T_{h_1} C^\ast_{\sigma_1}T_{g_1}^\ast  \qquad \mbox{and} \qquad C_{\psi}=T_{h_2} C^\ast_{\sigma_2}T_{g_2}^\ast,$$
where $g_1, h_1, \sigma_1$ and  $g_2, h_2, \sigma_2$ are  respectively the auxiliary functions for $\varphi$ and $\psi$ on Theorem A. Since
 $\varphi$ and $\psi$ are automorphisms  of  $D$, we have $\sigma_1=\varphi^{-1}$ and $\sigma_2=\psi^{-1}$.

For $w\in D$, let $K_w$ denote the reproducing kernel and $k_w$ be the normalized reproducing kernel respectively given by $$K_w(z)=\frac{1}{1-z\overline{w}} \qquad \mbox{and} \qquad k_w(z)=\frac{K_w(z)}{||K_w||}=\frac{\sqrt{1-|w|^2}}{1-z\overline{w}},$$ where $||\cdot||$ denotes the norm of  $H^2(D)$. Now, using the formula $C_\phi^{\ast}T^\ast_b K_w=\overline{b(w)}K_{\phi(w)}$, we get that{\setlength\arraycolsep{2pt}
\begin{eqnarray*} <C_\varphi K_z, C_\psi K_w> &=& <T_{h_1} C^\ast_{\sigma_1}T_{g_1}^\ast K_z, T_{h_2} C^\ast_{\sigma_2}T_{g_2}^\ast K_w>
\\ &=& g_1(z)\overline{g_2(w)}<T_{h_1} K_{\sigma_1(z)}, T_{h_2} K_{\sigma_2(w)}>
\\ &=& g_1(z)\overline{g_2(w)}< K_{\sigma_1(z)}, T^\ast_{h_1}T_{h_2} K_{\sigma_2(w)}>.
\end{eqnarray*}}

Since $h_1(z)=c_1z+ d_1$ and  $h_2(z)=c_2z+ d_2$ are in $ H^\infty$, using the semi-multiplicative property for Toeplitz operators mod $\mathcal{K}$ as mentioned in Section 2, we see
$$T^\ast_{h_1}T_{h_2}=T_{\overline{h_1}}T_{h_2}=T_{\overline{h_1}\, h_2}=T_{h_2\, \overline{h_1}}=T_{h_2}T_{\overline{h_1}}+L=T_{h_2}T^\ast_{h_1}+L,$$
where $L$ is a compact operator on $H^2(D)$. It follows that{\setlength\arraycolsep{2pt}
\begin{eqnarray*}&& <C_\varphi K_z, C_\psi K_w> =  g_1(z)\overline{g_2(w)}< K_{\sigma_1(z)}, T^\ast_{h_1}T_{h_2} K_{\sigma_2(w)}>
\\ &=&  g_1(z)\overline{g_2(w)}< K_{\sigma_1(z)}, T_{h_2}T^\ast_{h_1} K_{\sigma_2(w)}>+ g_1(z)\overline{g_2(w)}< K_{\sigma_1(z)}, L K_{\sigma_2(w)}>
\\ &=&  g_1(z)\overline{g_2(w)}< T_{h_2}^\ast K_{\sigma_1(z)},T^\ast_{h_1} K_{\sigma_2(w)}>+ g_1(z)\overline{g_2(w)}< K_{\sigma_1(z)}, L K_{\sigma_2(w)}>
\\ &=&  g_1(z)\overline{g_2(w)}h_2\circ\sigma_1(z) \overline{h_1\circ\sigma_2(w)}< K_{\sigma_1(z)}, K_{\sigma_2(w)}>
\\ &&+ g_1(z)\overline{g_2(w)}< K_{\sigma_1(z)}, L K_{\sigma_2(w)}>.
\end{eqnarray*}}

Fix $\omega\in \partial D$, since  $\varphi$ and $\psi$ are automorphisms, there exist $\zeta_1$ and $\zeta_2$ on $\partial D$ such that $\varphi(\zeta_1)=\psi(\zeta_2)=\omega$. Setting $z=r\zeta_2$ and $w=r\zeta_1$, then{\setlength\arraycolsep{2pt}
\begin{eqnarray*} &&\lim\limits_{r\to 1^-}<C_\varphi k_{r\zeta_2}, C_\psi k_{r\zeta_1}> = \lim\limits_{r\to 1^-}<C_\varphi \frac{K_{r\zeta_2}}{||K_{r\zeta_2}||}, C_\psi \frac{K_{r\zeta_1}}{||K_{r\zeta_1}||}>
\\ &=& \lim\limits_{r\to 1^-}(1-r^2)<C_\varphi K_{r\zeta_2}, C_\psi K_{r\zeta_1}>
\\ &=& \lim\limits_{r\to 1^-}(1-r^2) g_1(r\zeta_2)\,\overline{g_2(r\zeta_1)}\,h_2\circ\sigma_1(r\zeta_2)\, \overline{h_1\circ\sigma_2(r\zeta_1)}< K_{\sigma_1(r\zeta_2)}, K_{\sigma_2(r\zeta_1)}>
\\ && + \lim\limits_{r\to 1^-}(1-r^2)g_1(r\zeta_2)\,\overline{g_2(r\zeta_1)}< K_{\sigma_1(r\zeta_2)}, L K_{\sigma_2(r\zeta_1)}>
\\ &=&: I+II.
\end{eqnarray*}}Note that $\{k_w\}$ is a weakly convergent sequence  as $|w|\to 1$  with $||k_w||=1$ and $L$ is compact, we see that
$$\lim\limits_{|w|\to 1}\sqrt{1-|w|^2}||LK_w||=\lim\limits_{|w|\to 1}||Lk_w||=0.$$
Which gives that{\setlength\arraycolsep{2pt}
\begin{eqnarray*}&&\lim\limits_{r\to 1^-}\sqrt{1-r^2}||L K_{\sigma_2(r\zeta_1)}||
\\ &=& \lim\limits_{r\to 1^-}\sqrt{1-|\sigma_2(r\zeta_1)|^2}||L K_{\sigma_2(r\zeta_1)}||\cdot\sqrt{\frac{1-r^2}{1-|\sigma_2(r\zeta_1)|^2}}
\\ &=& 0\cdot\frac{1}{|\sigma_2'(\zeta_1)|^{1/2}}=0.
\end{eqnarray*}}
Hence,
{\setlength\arraycolsep{2pt}
\begin{eqnarray*} &&\lim\limits_{r\to 1^-}(1-r^2)|g_1(r\zeta_2)\,\overline{g_2(r\zeta_1)}< K_{\sigma_1(r\zeta_2)}, L K_{\sigma_2(r\zeta_1)}>|
\\ &\le & \lim\limits_{r\to 1^-}(1-r^2)|g_1(r\zeta_2)\,\overline{g_2(r\zeta_1)}|\cdot || K_{\sigma_1(r\zeta_2)}||\cdot  ||L K_{\sigma_2(r\zeta_1)}||
\\ &= & \lim\limits_{r\to 1^-}|g_1(r\zeta_2)\,\overline{g_2(r\zeta_1)}| \sqrt{\frac{1-r^2}{1-|\sigma_1(r\zeta_2)|^2}}\cdot \sqrt{1-r^2} ||L K_{\sigma_2(r\zeta_1)}||
\\ &= & |g_1(\zeta_2)\,\overline{g_2(\zeta_1)}|\frac{1}{|\sigma_1'(\zeta_2)|^{1/2}} \cdot 0=0
\end{eqnarray*}}and so $II=0$.

Now, using $\sigma_1=\varphi^{-1}$ and $\sigma_2=\psi^{-1}$, we calculate that{\setlength\arraycolsep{2pt}
\begin{eqnarray*} I &=& \lim\limits_{r\to 1^-}(1-r^2) g_1(r\zeta_2)\,\overline{g_2(r\zeta_1)}\,h_2\circ\sigma_1(r\zeta_2)\, \overline{h_1\circ\sigma_2(r\zeta_1)}< K_{\sigma_1(r\zeta_2)}, K_{\sigma_2(r\zeta_1)}>
\\ &=&
\lim\limits_{r\to 1^-} \overline{h_1\circ\psi^{-1}(r\zeta_1)}\ g_1(r\zeta_2)\ h_2\circ\varphi^{-1}(r\zeta_2)\ \overline{g_2(r\zeta_1)}\ \frac{1-r^2}{1-<\psi^{-1}(r\zeta_1), \varphi^{-1}(r\zeta_2)>}.
\end{eqnarray*}}It is easy to see that $$\lim\limits_{r\to 1^-}\frac{1-r^2}{1-<\psi^{-1}(r\zeta_1), \varphi^{-1}(r\zeta_2)>}=0$$ unless $\varphi^{-1}(\zeta_2)=\psi^{-1}(\zeta_1)$. Next, we assume $\varphi^{-1}(\zeta_2)=\psi^{-1}(\zeta_1)$. By Lemma 3.3,
 we get{\setlength\arraycolsep{2pt}
\begin{eqnarray*}\lim\limits_{r\to 1}\frac{1-r^2}{1-<\psi^{-1}(r\zeta_1), \varphi^{-1}(r\zeta_2)>}&=&\lim\limits_{r\to 1}<C_{\varphi^{-1}}^{\ast} k_{r\zeta_2}, C_{\psi^{-1}}^{\ast} k_{r\zeta_1}>
\\ &=&\frac{2}{|(\psi^{-1})'(\zeta_1)|+|(\varphi^{-1})'(\zeta_2)|}.
\end{eqnarray*}}On the other hand, an easy computation shows that
$$h_1(z)\overline{g_1\circ\varphi(z)}=\frac{c_1z+d_1}{\overline{a_1-c_1\varphi(z)}}
=\frac{c_1z+d_1}{\overline{a_1-c_1\frac{a_1z+b_1}{c_1z+d_1}}}=\frac{|c_1z+d_1|^2}{\overline{a_1d_1-c_1b_1}}=|c_1z+d_1|^2$$
and similarly
$$h_2(z)\overline{g_2\circ\psi(z)}=|c_2z+d_2|^2.$$
This fact has been proved by MacCluer and  Pons \cite{MP06} for the  the unit ball case. Combining this with $\varphi^{-1}(\zeta_2)=\psi^{-1}(\zeta_1)$, we obtain{\setlength\arraycolsep{2pt}
\begin{eqnarray*}\overline{h_1\circ\psi^{-1}(\zeta_1)}\ g_1(\zeta_2)&=& \overline{h_1\circ\varphi^{-1}(\zeta_2)}\ g_1\circ\varphi\circ\varphi^{-1}(\zeta_2)
\\ &=& |c_1\varphi^{-1}(\zeta_2)+d_1|^2= \biggl|c_1\frac{d_1\zeta_2-b_1}{-c_1\zeta_2+a_1}+d_1\biggr|^2
\\ &=& \frac{1}{|a_1-c_1\zeta_2|^2}=|(\varphi^{-1})'(\zeta_2)|
\end{eqnarray*}}
and{\setlength\arraycolsep{2pt}
\begin{eqnarray*}h_2\circ\varphi^{-1}(\zeta_2)\ \overline{g_2(\zeta_1)}&=&h_2\circ\psi^{-1}(\zeta_1)\ \overline{g_2\circ\psi\circ\psi^{-1}(\zeta_1)}
\\ &=& |c_2\psi^{-1}(\zeta_1)+d_2|^2=|(\psi^{-1})'(\zeta_1)|.
\end{eqnarray*}}As a result, we  show that $$I=|(\varphi^{-1})'(\zeta_2)||(\psi^{-1})'(\zeta_1)|\frac{2}{|(\psi^{-1})'(\zeta_1)|+|(\varphi^{-1})'(\zeta_2)|}$$
if $\varphi^{-1}(\zeta_2)=\psi^{-1}(\zeta_1)$; Otherwise $I=0$. Moreover,  since $\varphi(\zeta_1)=\psi(\zeta_2)=\omega\in\partial D$, if  $\varphi^{-1}(\zeta_2)=\psi^{-1}(\zeta_1)$ then we have $$\psi^{-1}\circ\varphi^{-1}(\omega)=\varphi^{-1}\circ\psi^{-1}(\omega).$$
This means $\varphi^{-1}$ and $\psi^{-1}$ commute and hence $\varphi$ and $\psi$ commute. Immediately, we get
$$(\psi^{-1})'(\varphi^{-1}(\omega))(\varphi^{-1})'(\omega)=(\varphi^{-1})'(\psi^{-1}(\omega))(\psi^{-1})'(\omega),$$
i.e. $$(\psi^{-1})'(\zeta_1)\frac{1}{\varphi'(\zeta_1)}=(\varphi^{-1})'(\zeta_2)\frac{1}{\psi'(\zeta_2)}.$$
Thus, the above discussions deduce that{\setlength\arraycolsep{2pt}
\begin{eqnarray*}&&\lim\limits_{r\to 1^-}<C_\varphi k_{r\zeta_2}, C_\psi k_{r\zeta_1}> =:I+II
\\ &=& |(\varphi^{-1})'(\zeta_2)| |(\psi^{-1})'(\zeta_1)|\frac{2}{|(\psi^{-1})'(\zeta_1)|+|(\varphi^{-1})'(\zeta_2)|}+0
\\ &=& \biggl|(\psi^{-1})'(\zeta_1)\frac{\psi'(\zeta_2)}{\varphi'(\zeta_1)}\biggr| |(\psi^{-1})'(\zeta_1)|\frac{2}{|(\psi^{-1})'(\zeta_1)|+ \biggl|(\psi^{-1})'(\zeta_1)\frac{\psi'(\zeta_2)}{\varphi'(\zeta_1)}\biggr| }
\\ &=&|(\psi^{-1})'(\zeta_1)\psi'(\zeta_2)| \frac{2}{|\varphi'(\zeta_1)|+|\psi'(\zeta_2)|}
\end{eqnarray*}}under the condition of $\varphi^{-1}(\zeta_2)=\psi^{-1}(\zeta_1)$. Otherwise, it is zero.

At last, by our hypothesis that $[C_\psi^{\ast}, C_\varphi]$ is compact on  $H^2(D)$, we see that{\setlength\arraycolsep{2pt}
\begin{eqnarray*}0 &=& \lim\limits_{r\to 1^-}|| [C_\psi^{\ast}, C_\varphi]k_{r\zeta_2}||\ge  \lim\limits_{r\to 1^-}|<[C_\psi^{\ast}, C_\varphi]k_{r\zeta_2},  k_{r\zeta_1}>|
\\ &=&
\lim\limits_{r\to 1^-}|<C_\varphi k_{r\zeta_2}, C_\psi k_{r\zeta_1}> -<C_{\psi}^{\ast} k_{r\zeta_2}, C_{\varphi}^{\ast} k_{r\zeta_1}>|.
\end{eqnarray*}}Using Lemma 3.3 again, we know that $$\lim\limits_{r\to 1^-}<C_{\psi}^{\ast} k_{r\zeta_2}, C_{\varphi}^{\ast} k_{r\zeta_1}>=\frac{2}{|\varphi'(\zeta_1)|+|\psi'(\zeta_2)|}.$$
It follows that $$\lim\limits_{r\to 1^-}<C_\varphi k_{r\zeta_2}, C_\psi k_{r\zeta_1}> =\lim\limits_{r\to 1^-}<C_{\psi}^{\ast} k_{r\zeta_2}, C_{\varphi}^{\ast} k_{r\zeta_1}>\ne 0.$$
Combining this with the previous conclusion, we must have  $\varphi^{-1}(\zeta_2)=\psi^{-1}(\zeta_1)$. Moreover, the above equality implies  $$|(\psi^{-1})'(\zeta_1)\psi'(\zeta_2)|=1.$$
(We don't know how to use this conclusion, but it maybe have some independent interest.) Note that we have obtained that  $\varphi$ and $\psi$ commute from
 $\varphi^{-1}(\zeta_2)=\psi^{-1}(\zeta_1)$. Therefore, $[C_\psi^{\ast}, C_\varphi]$ is compact must deduce that $\varphi$ and $\psi$ commute.
\ \ $\Box$
\\ \par In the proof of Lemma 3.4,  replacing  Lemma 3.3 by  Lemma 3.2,     we can  deduce  the similar result for the unit ball. We state it as the following and omit its proof.

\begin{lemma}\label{lem 3.5} If  $\varphi, \psi\in \mbox{Aut} (B_N)$  and  $[C_\psi^{\ast}, C_\varphi]$ is compact on  $H^2(B_N)$  or $A^2_s(B_N)$ ($s>-1$), then $\varphi$ and $\psi$ commute.
\end{lemma}

Next,  we will give a complete proof for  our main theorem. That is, we will show that when $\varphi$ and $\psi$ are automorphisms  of $B_N$, the compactness of   $[C_\psi^{\ast}, C_\varphi]$ implies that both $\varphi$ and $\psi$ are unitary. In the proof of this similar conclusion, the technique treated on  $A^2_s(D)$  (see \cite{MNW13}) is different from that on   $H^2(D)$ (see \cite{CLN12}). On  $A^2_s(D)$, the polar decomposition of $C_\varphi$ was used. I think which also holds on  the unit ball, so we
can use similar idea as on $A^2_s(D)$ to complete the proof of this result. However, in order to exhibit the special property of composition operator $C_\varphi$ when $\varphi$ is an automorphism of $B_N$, we try to use the following interesting  lemmas to prove this result.

\begin{lemma}\label{lem 3.6} Suppose that $\varphi$ is an automorphism of  $B_N$. Then on the space $\mathcal{H}$,
 $$C_{\varphi}^{\ast}=T_f C_{\varphi}^{-1}=T_f C_{\varphi^{-1}}, $$ where $T_f$ is the Toeplitz operator with symbol $$f(z)=\biggl(\frac{1-|\varphi(0)|^2}{|1-<z, \varphi(0)>|^2}\biggr)^t$$   with $t=N$ when  $\mathcal{H}=H^2(B_N)$  and $t=N+s+1$
when  $\mathcal{H}=A^2_s(B_N)$ ($s>-1$).
\end{lemma}

This result on the Hardy space  $H^2(B_N)$  and  the Bergman space  $A^2(B_N)$ was established by Bourdon and MacCluer \cite{BM07}. The extension  to  the weighted space $A^2_s(B_N)$  can be obtained similarly when using the change of variables formula in Proposition 1.13 of \cite{Zhu}. For the case of the unit disk, please see Theorem 4.2 of \cite{MNW13}.

\begin{lemma}\label{lem 3.7} Let $\varphi$ be an automorphism of  $B_N$. If $f$ is continuous on $\partial B_N$ or $\overline{B_N}$, then
 $$C_{\varphi}T_f- T_{f\circ \varphi} C_{\varphi} $$ is compact when acting on  $H^2(B_N)$ or   $A^2_s(B_N)$ ($s>-1$).
\end{lemma}

\proof Let $\varphi$ be an automorphism of  $B_N$ and $a=\varphi^{-1}(0)$. By Theorem 2.2.5 of \cite{Rudin}, the identity
$$1-<\varphi(z), \zeta>=1-<\varphi(z), \varphi\circ\varphi^{-1}(\zeta)>=\frac{(1-|a|^2)(1-<z, \varphi^{-1}(\zeta)>)}{(1-<z, a>)(1-<a, \varphi^{-1}(\zeta)>)}$$
holds for all $z\in \overline{B_N}$ and $\zeta\in \partial B_N$. Now, using this identity, for any $g\in H^2(B_N)$, we get that{\setlength\arraycolsep{2pt}
\begin{eqnarray*}&& (C_{\varphi}T_f \, g)(z)=(T_f g)(\varphi(z))=\int_{\partial B_N} f(\zeta) g(\zeta)\frac{1}{(1-<\varphi(z), \zeta>)^N}d\sigma(\zeta)
\\ &=& \int_{\partial B_N} f(\zeta) g (\zeta)\frac{(1-<z, a>)^N(1-<a, \varphi^{-1}(\zeta)>)^N}{(1-|a|^2)^N(1-<z, \varphi^{-1}(\zeta)>)^N}d\sigma(\zeta)
\\ &=& \int_{\partial B_N} f\circ\varphi(\eta) g\circ\varphi(\eta)\frac{(1-<z, a>)^N(1-<a, \eta>)^N}{(1-|a|^2)^N(1-<z, \eta>)^N}\cdot\frac{(1-|a|^2)^N}{|1-<\eta, a>|^{2N}}d\sigma(\eta)
\\ &=& (1-<z, a>)^N \int_{\partial B_N} f\circ\varphi(\eta) g\circ\varphi(\eta)\frac{1}{(1-<\eta, a>)^N}\cdot \frac{1}{(1-<z, \eta>)^N}d\sigma(\eta)
\\ &=& (1-<z, a>)^N (T_{f\circ\varphi\cdot  K_a}\,  g\circ\varphi)(z)
\\ &=& (T_{1/K_a}T_{f\circ\varphi} T_{ K_a} C_\varphi\, g)(z),
\end{eqnarray*}}where we have used the change of variables formula (see Corollary 4.4 of \cite{Zhu}) and the   kernel function
$$K_a(z)=\frac{1}{(1-<z, a>)^N}$$ and the function $1/K_a$ are  analytic on $\overline{B_N}$.  Hence, $$C_{\varphi}T_f =T_{1/K_a}T_{f\circ\varphi} T_{ K_a} C_\varphi.$$

Since $f$ is continuous on $\partial B_N$, by the semi-multiplicative property for Toeplitz operators mod $\mathcal{K}$ in Section 2, we know that
$$T_{1/K_a}T_{f\circ\varphi}= T_{1/K_a\cdot f\circ\varphi}+L=T_{f\circ\varphi \cdot 1/K_a}+L= T_{f\circ\varphi} T_{1/K_a}+L,$$
where $L$ is a compact operator on $H^2(B_N)$. Therefore,
 $$C_{\varphi}T_f =T_{1/K_a}T_{f\circ\varphi} T_{ K_a} C_\varphi=(T_{f\circ\varphi} T_{1/K_a}+L)T_{ K_a} C_\varphi=T_{f\circ\varphi} C_\varphi+L', $$
where $L'$ is  compact. Applying similar technique, we can obtain the similar result for the weighted Bergman space and so we complete the proof.
\ \ $\Box$
\\ \par
Based on these lemmas, we will use another technique to prove Theorem 3.1, which generalizes Theorem 5.2 of \cite{CLN12} and Theorem 5.1 of \cite{MNW13} to the unit ball.
\\ \\
{\bf Proof of Theorem 3.1.} We only need to prove the "only if" part and we only give a proof for the  Hardy space case.

Assume that   $[C_\psi^{\ast}, C_\varphi]$ is compact. Since  $\psi$ is an automorphism of  $B_N$, by Lemma 3.6, we have $$C_{\psi}^{\ast}=T_f C_{\psi^{-1}} $$ with  $$f(z)=\biggl(\frac{1-|\psi(0)|^2}{|1-<z, \psi(0)>|^2}\biggr)^N.$$
Thus, $$[C_\psi^{\ast}, C_\varphi]=C_\psi^{\ast} C_\varphi- C_\varphi C_\psi^{\ast}=T_f C_{\psi^{-1}}C_\varphi- C_\varphi T_f C_{\psi^{-1}}.$$
Now, using Lemma 3.5,   $[C_\psi^{\ast}, C_\varphi]$ is compact means that $\varphi$ and  $\psi$ commute, i.e.  $$\varphi\circ\psi=\psi\circ\varphi.$$
This gives $\varphi=\psi\circ\varphi\circ\psi^{-1}$ and{\setlength\arraycolsep{2pt}
\begin{eqnarray*}[C_\psi^{\ast}, C_\varphi]C_\psi &=&(T_f C_{\psi^{-1}}C_\varphi- C_\varphi T_f C_{\psi^{-1}})C_\psi
\\ &=& T_f C_{\psi^{-1}}C_\varphi C_\psi- C_\varphi T_f C_{\psi^{-1}} C_\psi
\\ &=& T_f C_{\psi\circ\varphi\circ\psi^{-1}}- C_\varphi T_f
\\ &=& T_f C_{\varphi}- C_\varphi T_f.
\end{eqnarray*}}

Note that $f$ is continuous on $\overline{B_N}$ and $\varphi$ is an automorphism of  $B_N$. It follows from Lemma 3.7 that
$$C_{\varphi}T_f\equiv T_{f\circ \varphi} C_{\varphi} \quad (\mbox{mod} \ \mathcal{K}).$$ Therefore,  {\setlength\arraycolsep{2pt}
\begin{eqnarray*}[C_\psi^{\ast}, C_\varphi]C_\psi &=& T_f C_{\varphi}- C_\varphi T_f
\\ &\equiv & T_f C_{\varphi}- T_{f\circ \varphi} C_{\varphi}\quad (\mbox{mod} \ \mathcal{K})
\\ &= & T_{f - f\circ \varphi} C_{\varphi}\quad (\mbox{mod}\  \mathcal{K}).
\end{eqnarray*}}Finally, since $C_\psi$ and  $C_\varphi$  are invertible, we see that $[C_\psi^{\ast}, C_\varphi]$ is compact if and only if  $[C_\psi^{\ast}, C_\varphi]C_\psi$ is compact, which is equivalent to that  $ T_{f - f\circ \varphi}$ is compact. Applying Lemma 2 of  \cite{Cob73}, we see that $f - f\circ \varphi\equiv 0$ on  $\overline{B_N}$. Since $\varphi$ is not the identity, the representation of $f$ gives that $f$ must be constant. Thus, we get  $\psi(0)=0$.

As we know, an automorphism $\phi$ of $B_N$ is a unitary transformation of $\mathbb{C}^N$ if and only if $\phi(0)=0$ (see Lemma 1.1 of \cite{Zhu}). Together this with the above discussion, we see that $\psi$ must be unitary. On the other hand, $[C_\psi^{\ast}, C_\varphi]=C_\psi^{\ast} C_\varphi- C_\varphi C_\psi^{\ast}$
is compact implies that $$(C_\psi^{\ast} C_\varphi- C_\varphi C_\psi^{\ast})^\ast=C_\varphi^{\ast} C_\psi- C_\psi C_\varphi^{\ast}=[C_\varphi^{\ast}, C_\psi]$$
is compact. So similar arguments  deduce that $\varphi$ is also unitary. Consequently, all these give that both $\varphi$ and  $\psi$ are  unitary and they commute under the condition that $[C_\psi^{\ast}, C_\varphi]$ is compact. \ \ $\Box$


\section {The commutator on the Dirichlet space}

 In this section, we try to characterize the compactness of $[C_\psi^{\ast}, C_\varphi]$  on the Dirichlet space $\mathcal{D}(B_N)$, where $\varphi$ and $\psi$ are linear fractional self-maps  of $B_N$. The following lemma is about compact difference of  linear fractional composition operators on
$\mathcal{D}(B_N)$, which will prove useful for our result.

\begin{lemma}\label{lem 4.1} Suppose  that $\varphi$ and $\psi$ are linear fractional self-maps of $B_N$. Then $C_\varphi-C_\psi$ is compact on $\mathcal{D}(B_N)$ if and only if $\varphi=\psi$ or both $C_\varphi$ and $C_\psi$ are compact.
\end{lemma}

\proof This result for a special case has been pointed out by Pons \cite{Po11}, without proof. For completeness, we give a simple proof.

In \cite{Po11}, for real $s$, the weighted Dirichlet space $\mathcal{D}_s(B_N)$ is defined by
$$\mathcal{D}_s(B_N)=\{f(z)=\sum\limits_\alpha c_\alpha z^\alpha\ \mbox{analytic in}\ B_N: \sum\limits_\alpha(|\alpha|+1)^{1-s}|c_\alpha|^2\omega_\alpha<\infty\},$$
where $$\omega_\alpha=||z^\alpha||^2=\frac{(N-1)!\alpha!}{(N-1+|\alpha|)!}.$$
On the other hand, let $\beta(k)=(k+1)^t$ for real $t$, if $f(z)=\sum\limits_\alpha c_\alpha z^\alpha=\sum^\infty_0f_k(z)$ is analytic in $B_N$, then $f$ belongs to the weighted Hardy space $H^2(\beta, B_N)$ if and only if $$\sum^\infty_0||f_k||^2 \beta(k)^2=\sum\limits_\alpha(|\alpha|+1)^{2t}|c_\alpha|^2\omega_\alpha<\infty.$$
Thus, the weighted Dirichlet space $\mathcal{D}_s(B_N)$, in fact, is the weighted Hardy space $H^2(\beta, B_N)$ with the weight $\beta(k)=(k+1)^{(1-s)/2}$. Using Proposition 2.4 of \cite{JC14}, we see that all linear fractional self-maps induce bounded composition operators on  $\mathcal{D}_s(B_N)$.

Suppose $s_1<s<s_2<\infty$, complex interpolation theorem for the  weighted Dirichlet space $\mathcal{D}_s(B_N)$ tells us $$[\mathcal{D}_{s_1}, \mathcal{D}_{s_2}]_\theta=\mathcal{D}_s$$ with $s=(1-\theta)s_1+\theta s_2$
for $\theta\in (0, 1)$ (see Proposition 1 of \cite{Po11}). Choosing $s_1=-n$ and $s_2=2$, then  for linear fractional self-maps $\varphi$ and $\psi$ of $B_N$,
the operator $C_\varphi-C_\psi$
is bounded on $\mathcal{D}_{-n}(B_N)$ and $\mathcal{D}_2(B_N)=A^2(B_N)$. Now, if $C_\varphi-C_\psi$ is compact on $\mathcal{D}(B_N)=\mathcal{D}_{1-n}(B_N)$. Since $1-n=(1-\theta)(-n)+\theta\cdot 2$ with $\theta=\frac{1}{n+2}$,  using the compactness theorem for interpolating operators (see Theorem 2.1 in \cite{Cwi92}), then $C_\varphi-C_\psi$ is also compact on $H^2(B_N)=\mathcal{D}_1(B_N)$ with $1=(1-\theta)(-n)+\theta\cdot 2$ for $\theta=\frac{n+1}{n+2}$. Here, all above spaces are identified to equal with an equivalent norm. Applying Theorem 2 of \cite{HMW11} or Theorem 3.1 of \cite{JO11}, we get that $\varphi=\psi$ or both $C_\varphi$ and $C_\psi$ are compact. So we have shown one direction. Another direction is obvious and  the proof is completed. \ \ $\Box$
\\ \par
When $\varphi$ is a linear fractional self-map of $B_N$, the adjoint of  composition operator $C_\varphi$ on $\mathcal{D}(B_N)$ is  mainly determined by another composition operator. Together this adjoint property with Lemma 4.1, we will give the following condition for  $[C_\psi^{\ast}, C_\varphi]$ to be compact on  $\mathcal{D}(B_N)$.

\begin{theorem} \label{the 4.2} Let $\varphi$  and $\psi$ be  linear fractional self-maps of  $B_N$. If  $[C_\psi^{\ast}, C_\varphi]\ne 0$,  then $[C_\psi^{\ast}, C_\varphi]$ is non-trivially  compact on  $\mathcal{D}(B_N)$ if and only if   $||\psi||_\infty=||\varphi||_\infty=1$  and $\varphi\circ\sigma=\sigma\circ\varphi$,  where $\sigma$ is the Krein adjoint  of $\psi$.
\end{theorem}

\proof For $w\in B_N$, let $K_w$ denote the reproducing kernel of $\mathcal{D}(B_N)$. Theorem C gives   $$C_{\psi}^{\ast}f=f(0)K_{\psi(0)}+C_\sigma f-f(\sigma(0))$$ for any $f\in \mathcal{D}(B_N)$, where $\sigma$ is the Krein adjoint  of $\psi$.
Thus, {\setlength\arraycolsep{2pt}
\begin{eqnarray*}[C_\psi^{\ast}, C_\varphi]f &=& (C_\psi^{\ast} C_\varphi- C_\varphi C_\psi^{\ast}) f=C_\psi^{\ast} C_\varphi f- C_\varphi C_\psi^{\ast} f
\\ &=& f\circ\varphi(0) K_{\psi(0)}-C_\sigma  C_\varphi f-f\circ\varphi\circ\sigma(0)
\\ && - [f(0)K_{\psi(0)}\circ\varphi- C_\varphi C_\sigma f-f(\sigma(0))]
\\ &=& (C_\varphi C_\sigma-C_\sigma  C_\varphi )f+f\circ\varphi(0) K_{\psi(0)}+f(\sigma(0))
\\ && -f\circ\varphi\circ\sigma(0)-f(0)K_{\psi(0)}\circ\varphi.
\end{eqnarray*}}This implies that $[C_\psi^{\ast}, C_\varphi]$ is compact  on  $\mathcal{D}(B_N)$  if and only if $$[C_\varphi, C_\sigma]=C_\varphi C_\sigma-C_\sigma  C_\varphi=C_{\sigma\circ\varphi}-C_{\varphi\circ\sigma}$$ is compact.

It is easy to see that $C_\psi^{\ast} C_\varphi$ is compact if and only if $C_\sigma  C_\varphi =C_{\varphi\circ\sigma}$ is compact, and the compactness of $C_\varphi C_\psi^{\ast}$ is equivalent to the compactness of $C_\varphi C_\sigma=C_{\sigma\circ\varphi}$. Since $[C_\psi^{\ast}, C_\varphi]\ne 0$,  we get  that  $[C_\psi^{\ast}, C_\varphi]$ is non-trivially  compact  on  $\mathcal{D}(B_N)$ if and only if
$C_{\sigma\circ\varphi}-C_{\varphi\circ\sigma}$ is also  non-trivially  compact. By Lemma 4.1, this is equivalent to $\varphi\circ\sigma=\sigma\circ\varphi$ and $||\varphi\circ\sigma||_\infty=||\sigma\circ\varphi||_\infty=1$, i.e.   $||\psi||_\infty=||\varphi||_\infty=1$. So we  complete the proof.  \ \ $\Box$
\\ \par Note that if $\psi$ is an automorphism of $B_N$, then $\sigma=\psi^{-1}$. Thus, $\varphi\circ\sigma=\sigma\circ\varphi$ gives that $\varphi\circ\psi^{-1}=\psi^{-1}\circ\varphi$. This is the same to $\varphi\circ\psi=\psi\circ\varphi$. Immediately, as a corollary of Theorem 4.2, we obtain the following result.

\begin{theorem} \label{the 4.3} If $\varphi$  and $\psi$  are automorphisms of  $B_N$. Then $[C_\psi^{\ast}, C_\varphi]$ is   compact on  $\mathcal{D}(B_N)$ if and only if $\varphi$ and $\psi$  commute.
\end{theorem}

Using Theorem 4.2  and similar discussions in the proof of Theorem 3.1 in \cite{MNW13}, we also obtain the following similar result on the Dirichlet space $\mathcal{D}(D)$
as on the Hardy space $H^2(D)$ and the weighted Bergman space $A^2_s(D)$ ($s>-1$).

\begin{theorem} \label{the 4.4} Suppose that $\varphi$  and $\psi$ are   linear fractional self-maps of  $D$, and one of which is a non-automorphism.  The commutator $[C_\psi^{\ast}, C_\varphi]$ is non-trivially  compact  on  $\mathcal{D}(D)$ if and only if one of the following is true.
\\ (i)\  $\varphi$  and $\psi$ are both parabolic with the same boundary fixed point,
\\ (ii)\ $\varphi$  and $\psi$ are hyperbolic such that the fixed points of $\varphi$  and $\psi$ are $(\zeta, a)$ and $(\zeta, 1/\overline{a})$
respectively
with $\zeta\in \partial D$.
\end{theorem}

\noindent {\bf Remark.} We first see the following example. Let
{\setlength\arraycolsep{2pt}
\begin{eqnarray*} \varphi(z)=\left(
\begin{array}{cc}
 1 & 0
\\ 0  & 1/2
\end{array}\right)\left(
\begin{array}{c}
z_1
\\ z_2
\end{array}\right)=(z_1, z_2/2), \qquad z=(z_1, z_2)\in B_2
\end{eqnarray*}}and{\setlength\arraycolsep{2pt}
\begin{eqnarray*} \psi(z)=\left(
\begin{array}{cc}
 1 & 0
\\ 0  & 1/3
\end{array}\right)\left(
\begin{array}{c}
z_1
\\ z_2
\end{array}\right)=(z_1, z_2/3), \qquad z=(z_1, z_2)\in B_2.
\end{eqnarray*}}Thus, {\setlength\arraycolsep{2pt}
\begin{eqnarray*} \sigma(z)=\left(
\begin{array}{cc}
 1 & 0
\\ 0  & 1/3
\end{array}\right)^\ast\left(
\begin{array}{c}
z_1
\\ z_2
\end{array}\right)=\left(
\begin{array}{cc}
 1 & 0
\\ 0  & 1/3
\end{array}\right)\left(
\begin{array}{c}
z_1
\\ z_2
\end{array}\right)
\end{eqnarray*}}and so{\setlength\arraycolsep{2pt}
\begin{eqnarray*} \varphi\circ\sigma(z)=\left(
\begin{array}{cc}
 1 & 0
\\ 0  & 1/6
\end{array}\right)\left(
\begin{array}{c}
z_1
\\ z_2
\end{array}\right)=\sigma\circ\varphi(z).
\end{eqnarray*}}Moreover, it is easy to check that{\setlength\arraycolsep{2pt}
\begin{eqnarray*}[C_\psi^{\ast}, C_\varphi]f &=&  (C_\varphi C_\sigma-C_\sigma  C_\varphi )f+f\circ\varphi(0) K_{\psi(0)}+f(\sigma(0))
\\ && -f\circ\varphi\circ\sigma(0)-f(0)K_{\psi(0)}\circ\varphi
\\ &=& 0
\end{eqnarray*}}for any $f\in \mathcal{D}(B_2)$.
As a consequence, when $\varphi\circ\sigma=\sigma\circ\varphi$ and  $||\psi||_\infty=||\varphi||_\infty=1$, it may happen that $[C_\psi^{\ast}, C_\varphi]=0$. Hence, in  Theorem 4.2, we need assume that  $[C_\psi^{\ast}, C_\varphi]\ne 0$. However, from the proof of Theorem 3.1 in \cite{MNW13}, this condition is not necessary
for the unit disk.
\\ \\ \\
\centerline{ ACKNOWLEDGEMENTS }
\\ \par During the author's visiting at The College at Brockport, State University of New York, they provided a good environment for working on this paper. Shanghai Municipal Education Commission  provided the financial support during her visiting. She would like to express her gratitude to them.

\bibliographystyle{amsplain}

\end{document}